\documentclass[12pt]{amsart}

\usepackage{amssymb,amsthm,amsmath,verbatim}
\usepackage{hyperref}
\usepackage[alphabetic]{amsrefs}
\usepackage[all]{xy}
\usepackage{enumerate}

\usepackage[margin=3cm]{geometry}

\newtheorem{theorem}{Theorem}[section]
\newtheorem{corollary}[theorem]{Corollary}
\newtheorem{lemma}[theorem]{Lemma}
\newtheorem{proposition}[theorem]{Proposition}
\newtheorem{definition}[theorem]{Definition}

\theoremstyle{remark}

\newtheorem{question}[theorem]{Question}

\newcommand{\N}{\mathbb{N}}

\newcommand{\C}{\mathbb{C}}

\newcommand{\cuntzeq}{\sim}
\newcommand{\cuntzle}{\preccurlyeq}
\newcommand{\Cu}{\mathrm{Cu}}

\newcommand{\QT}{\mathrm {QT}}

\newcommand{\ZZ}{\mathcal{Z}}

\renewcommand{\epsilon}{\varepsilon}
\renewcommand{\leq}{\leqslant}
\renewcommand{\geq}{\geqslant}

\title{Remarks on $\ZZ$-stable projectionless C*-algebras}
\author{Leonel Robert}

\begin{document}
\begin{abstract}
It is shown that $\ZZ$-stable  projectionless C*-algebras have the property that
every element is a limit of products of two nilpotents.
This is then used  to classify the approximate unitary equivalence classes of positive elements in 
such C*-algebras using traces. 
\end{abstract}

\maketitle


\section{Introduction}
Let us denote by $\ZZ$ the Jiang-Su C*-algebra (\cite{jiang-su}). A C*-algebra $A$ is called $\ZZ$-stable if $A\cong A\otimes \ZZ$.
A theme that has unfolded in the past decade in the field of ``structure and classification of C*-algebras"   is how the $\ZZ$-stable C*-algebras have good properties   that separate them from the ``pathological" examples found by Villadsen, R\o rdam, Toms  and others (see \cite{rordam, brown-perera-toms, elliott-toms}).
The results of this note contribute to the study of $\ZZ$-stable projectionless C*-algebras.   Henceforth, the term ``projectionless" is used to designate the C*-algebras none of whose quotients contains a non-zero projection. We prove the following theorem.

\begin{theorem}\label{nilpotentsTH}
Let $A$ be a $\ZZ$-stable projectionless C*-algebra. Then for every   $x\in A$ there exist nilpotent
elements $y_n,z_n$, with $n=1,2,\dots$, in the hereditary subalgebra generated by $x$, such that $y_nz_n\to x$.
\end{theorem}
Theorem \ref{nilpotentsTH} implies that every element of a $\ZZ$-stable projectionless C*-algebra can be approximated by invertible elements of its unitization. This ``almost stable rank one"  property has several interesting consequences. It implies, for example, the equivalence of (i) and (ii) in the following theorem:
\begin{theorem}\label{cuntzequiv}
Let $A$ be a $\ZZ$-stable projectionless C*-algebra. Let $a,b\in A_+$. The following are equivalent:
\begin{enumerate}[(i)]
\item
$a$ is Cuntz smaller than $b$,
\item
$\overline{aA}\hookrightarrow \overline{bA}$  as right $A$-Hilbert C*-modules,
\item
$d_\tau(a)\leq  d_\tau(b)$ for all lower semicontinuous 2-quasitraces $\tau\colon A_+\to [0,\infty]$.
\end{enumerate}
Furthermore, if $d_\tau(a)= d_\tau(b)$ for all lower semicontinuous 2-quasitraces $\tau$ then $\overline{aA}\cong \overline{bA}$
as right $A$-Hilbert C*-modules.
\end{theorem}
The hypothesis of $\ZZ$-stability cannot be dropped in Theorem \ref{cuntzequiv}, even in the simple nuclear case, as demonstrated in \cite{tikuisis}. In \cite{nawata1}, Nawata uses Theorem \ref{cuntzequiv} in the investigation of Picard groups of $\ZZ$-stable  projectionless simple C*-algebras.

The classification result from \cite{robert-santiago} is also applicable under the almost stable rank one property. We thus deduce the following theorem:
\begin{theorem}\label{positiveorbits}
Let $A$ be a $\ZZ$-stable projectionless C*-algebra. Let $a,b\in A_+$. The following are equivalent:
\begin{enumerate}[(i)]
\item
$a$ is approximately unitarily equivalent to $b$,

\item $\tau(f(a))=\tau(f(b))$ for any lower semicontinuous 2-quasitrace $\tau\colon A_+\to [0,\infty]$ and all $f\in C_0(0,1]_+$.

\item
$d_\tau((a-t)_+)=d_\tau((b-t)_+)$ for all  $t\geq 0$ and all $\tau$ as in (ii).
\end{enumerate}
\end{theorem}

\emph{Acknowledgement.} I am grateful to Bhishan Jacelon for fruitful conversations on the 
subject of this paper.

\section{Proof of Theorem \ref{nilpotentsTH}}
Let us first introduce some notation.
Let $A$ be a C*-algebra. Let us denote by $A_+$ the positive elements of $A$.
We will make use of the Cuntz comparison relations on $A_+$. If $a,b\in A_+$ then we write $a\precsim b$ if $d_n^*bd_n\to a$ for some $d_n\in A$
and  $a\sim b$ if  $a\precsim b$ and $b\precsim a$. 

Let $\mathcal K$ denote the C*-algebra of compact operators on a separable Hilbert space.
The Cuntz semigroup $\Cu(A)$
is defined as the set $(A\otimes \mathcal K)_+/\sim$ endowed with a suitable order and addition operation (see \cite{ara-perera-toms}).  Given $a\in (A\otimes \mathcal K)_+$, we shall denote by $[a]\in \Cu(A)$ its Cuntz equivalence class.

\begin{lemma}\label{keylemma}
Let $A$ be a   C*-algebra and let $x\in A$ be such that there exist $e,f\in A_+$ with the properties that $ex=xe=x$, $ef=0$, and $f$ is a full element of $A$. Then $x\otimes 1\in A\otimes \ZZ$ is the product
of two nilpotent elements. 
\end{lemma}

\begin{proof}
The relations $ex=xe=x$ imply that $x$ belongs to the Pedersen ideal of $A$. Therefore, since $f$ is full, there exist $a_i,b_i\in A$, with $i=1,2,\dots,m$, such that $x=\sum_{i=1}^m a_ifb_i$. Multiplying by $e$ on the left and on the right if necessary, we may assume that $fa_i=b_if=0$ for all $i$.

Let $n\geq 2m$. It is possible to find positive elements $e_1,e_2,\dots,e_n,d\in \ZZ$ such that 
\begin{enumerate}
\item
$\sum_{j=1}^n e_j+d=1$, 
\item
the elements $e_1,e_2,\cdots,e_n$ are pairwise orthogonal,
\item
there exist $(w_j)_{j=1}^{n-1}$ in $\ZZ$ such that $e_j=w_jw_j^*$ and $e_{j+1}=w_jw_j^*$ for all $j$,
\item
$[d]\leq [e_1]$ in the Cuntz semigroup $\Cu(\ZZ)$ of $\ZZ$. 
\end{enumerate}
(In fact, these elements can be found in the dimensionm drop algebra $Z_{n-1,n}$, which embeds in $\ZZ$; see \cite[Lemma 4.2]{rordam}).   Since $m[e_{n}+d]\leq 2m[e_1]\leq [1]$ in $\Cu(\ZZ)$, there exist $v_i\in \ZZ$, with $i=1,\dots,m$, such that $e_n+d=v_iv_i^*$ for all $i$ and the elements $v_i^*v_i\in \ZZ$ are pairwise orthogonal for all $i$. Let us now define $\alpha,\beta,\gamma,\delta\in A\otimes \ZZ$ 
as follows:
\begin{alignat*}{2}
\alpha &=  \sum_{i=1}^m a_jf^{1/2}\otimes v_j,\quad& \beta &=\sum_{i=1}^m f^{1/2}b_j\otimes v_j^*,\\
\gamma &= \sum_{j=1}^{n-1} x\otimes w_j,\quad 
&\delta &=\sum_{j=1}^{n-1} e\otimes w_j^*.
\end{alignat*}

We have $\gamma\beta=\alpha\delta=0$ (by the orthogonality of $e$ and $f$). Therefore,
\begin{align*}
(\gamma+\alpha)(\delta+\beta) &=\gamma\delta+\alpha\beta\\
&=\sum_{j=1}^{n-1} x\otimes e_i+(\sum_{i=1}^m a_ifb_i)\otimes (e_n+d)\\
&=x\otimes 1.
\end{align*} 

Let us now show that $\gamma+\alpha$ and $\delta+\beta$ are nilpotent elements.
We have that $\gamma^n=0$ and $\alpha^2=0$ (since $fa_i=0$ for all $i$), Finally, we have that
$\alpha\gamma=0$. Hence,
$
(\gamma+\alpha)^k=\sum_{i=1}^k\gamma^i\alpha^{k-i} 
$
for all $k$. Thus, for $k=n+1$ we get $(\gamma+\alpha)^{n+1}=0$. Similarly, $\delta$ and $\beta$ are nilpotent and $\delta\beta=0$. Thus,  $\delta+\beta$ is nilpotent.
\end{proof}

\begin{proof}[Proof of Theorem \ref{nilpotentsTH}]
Let us  identify $A$ with $A\otimes \ZZ$. 
Every element of $A\otimes\ZZ$ is approximately unitarily equivalent to one of the form $x\otimes 1$ (see the proof of \cite[Theorem 5.5]{brown-perera-toms}).
Thus, it suffices to assume that the given element has the form $x\otimes 1$. Set $x^*x+xx^*=a$. 
Since the property of being projectionless passes to hereditary subalgebras, we may assume that $a$ generates $A$ as a hereditary suablgebra (i.e., $a$ is strictly positive).
Let us choose an approximate unit $e_n\in C^*(a)$ of $A$ such that $e_{n+1}e_n=e_n$ for all $n\in \N$. 
Let us set $e_nxe_n=x_n$. Since $x_n\to x$, it suffices to show  that  $x_n\otimes 1$ is the product of two nilpotents for all $n\in \N$. Let $f_n\in C^*(a)_+$ be such that $e_nf_n=0$ and  $\delta_n a\leq e_{n+1}+f_n$ for some scalar $\delta_n>0$ and all $n$. The desired conclusion will follow from the previous lemma once we have shown that $f_n$ is a full element for all $n\in \N$. Let us fix $n\in \N$ and let $I$ denote the closed two-sided ideal generated by $f_n$. Let us suppose for the sake of contradiction that  $I\neq A$. In $A/I$, we have that $\delta_n \overline{a}\leq \overline{e}_{n+1}$, where $\overline{a},\overline{e}_{n+1}\in A/I$ denote the images of $a$ and $e_{n+1}$. This implies that $0$ is an isolated point of the spectrum of $\overline{a}$, which in turn implies that $A/I$ contains a non-zero projection. But this contradicts the assumption that $A$ is projectionless. Thus,
$f_n$ is full for all $n\in \N$.
\end{proof}

\section{Proofs of  Theorems \ref{cuntzequiv} and \ref{positiveorbits}}
\begin{definition}
Let $A$ be a C*-algebra. Let us say that $A$ almost has stable rank one if for every
 closed hereditary subalgebra $B\subseteq A$ we have $B\subseteq \overline{\mathrm{GL}(B^\sim)}$.
\end{definition}

\begin{corollary}\label{almostsr1}
Let $A$ be a $\ZZ$-stable projectionless C*-algebra. Then $A$ almost has stable rank one.
\end{corollary}
\begin{proof}
This follows from Theorem \ref{nilpotentsTH} and the fact that $x+\lambda\cdot  1$ is invertible
if $x$ is nilpotent and $\lambda\in \C\backslash\{0\}$.
\end{proof}

Before proving Theorems \ref{cuntzequiv} and Theorem \ref{positiveorbits},  let us introduce notation and make some preparatory remarks.

Let $A$ be a C*-algebra.
Let  us denote by $\QT_2(A)$  the cone of lower semicontinuous 2-quasitraces on 
the C*-algebra $A$. A lower semicontinuous 2-quasitrace
$\tau\in \QT_2(A)$ induces a dimension function $d_\tau$ on the positive elements of $A$
given by $d_\tau(a)=\lim_n \tau(a^{1/n})$ for $a\in A^+$. The value of $d_\tau$ on $a$ depends only on the Cuntz class
of $a$.

Let  $a,b\in A_+$.
If $A$ has stable rank one, then $a\precsim b$ if and only if  $\overline{aA}$ embeds as a Hilbert C*-module over $A$ in  $\overline{bA}$ and $a\sim b$ if and only if  $\overline{aA}$ is isomorphic to  $\overline{bA}$. This was first shown in \cite[Theorem 3]{coward-elliott-ivanescu}
using the language of Hilbert C*-modules. It was  later re-proven using positive elements in \cite[Proposition 1]{ciuperca-elliott-santiago} and \cite[Proposition 1.5]{huaxinlin}. By inspecting the  proofs in \cite{coward-elliott-ivanescu} and 
\cite{ciuperca-elliott-santiago}, it can be seen that they rely only
on the fact that $B\subseteq \overline{\mathrm{GL}(B^\sim)}$ for every $\sigma$-unital hereditary suablgebra 
$B$ of $A$. Thus, we arrive at the following proposition:

\begin{proposition}\label{Cuntzisomorphism}
Let $A$ be a C*-algebra. Let $a,b\in A_+$. Suppose that $A$  almost has stable rank one. 
Then $a\cuntzle b$ if and only if
$\overline{aA}\hookrightarrow \overline{bA}$ and $a\cuntzeq b$ if and only if $\overline{aA}\cong  \overline{bA}$.
\end{proposition}

Let $a\in A_+$.
The Cuntz semigroup  element $[a]$ is called compact if it is compactly contained in itself; i.e., $[a]\ll [a]$. It is shown in \cite{brown-ciuperca} that if  $[a]$ is compact then either $0$ is an isolated point in the spectrum of $a$ or $A$ contains  a scaling element (see \cite[Lemma 3.1 and Proposition 3.2]{brown-ciuperca}). In either case, $A$ contains a non-zero projection. It follows that if $A$ is projectionless and $a\in A_+$ then neither $[a]$ not its image after passing to a quotient of $A$ can be  a non-zero compact element. Thus, $[a]$ is a purely non-compact element in the sense of \cite{ers}.

\begin{proof}[Proof of Theorem \ref{cuntzequiv}]
The equivalence if (i) and (ii)  and the fact that Cuntz equivalence implies isomorphism of the right ideals generated by $a$ and $b$,  follow from Corollary \ref{almostsr1} and Propoposition \ref{Cuntzisomorphism}. 
The implication (i)$\Rightarrow $(iii) is well known. Finally,  let us show (iii)$\Rightarrow$(i). As it was argued in the previous paragraph, every non-zero positive element of $A$ gives rise to a purely non-compact element of  $\Cu(A)$. But it is shown in \cite[Theorem 6.6]{ers} that if $A$ is $\ZZ$-stable and $[a]$ and $[b]$ are purely non-compact elements such that $d_\tau(a)\leq d_\tau(b)$ for any $\tau\in \QT_2(A)$ then $[a]\leq [b]$.
This concludes the proof.
\end{proof}

\begin{proof}[Proof of Theorem \ref{positiveorbits}]
The implications (i)$\Rightarrow$(ii)$\Rightarrow$(iii) are clear. Let us prove  (iii)$\Rightarrow$(i).
By the previous proposition, $[(a-t)_+]=[(b-t)_+]$ for all $t\geq 0$.  By the classification result \cite[Theorem 1]{robert-santiago}, (i) will follow once we have shown that for every $x,e\in A$, with $e$ a positive contraction such that $ex=xe=x$, the elements $x^*x+e$ and $xx^*+e$ are stably approximately unitarily equivalent. 
It is shown in \cite[Proposition 4 (i)]{robert-santiago} that C*-algebras of stable rank one have this property. Exactly the same proof applies to the case that $A$ is $\ZZ$-stable and projectionless, since all that the proof uses is the almost stable rank one property defined above. 
\end{proof}

It is not true that a $\ZZ$-stable projectionless C*-algebra must have stable rank one.  It is worth noting, however, that in \cite{santiago} Santiago shows that tensoring an approximately subhomogeneous  C*-algebra by the Jacelon-Razak algebra  results in a ($\ZZ$-stable, projectionless) C*-algebra of  stable rank one.
The following question remains open:
\begin{question}
Let $A$ be simple, stably projectionless, and $\ZZ$-stable. Is $A$ of stable rank one?
\end{question}
In view of Rordam's \cite[Theorem 6.7]{rordam} asserting that a finite simple unital $\ZZ$-stable C*-algebra has stable rank one,  the answer to the above question is most likely positive.

\end{document}